\documentclass[12pt, a4paper]{article}
\usepackage[latin1]{inputenc}
\usepackage{vmargin}
\usepackage{amsfonts}
\usepackage{amssymb}
\usepackage{epstopdf}
\usepackage[dvips]{graphicx}
\usepackage{amsmath}
\usepackage{amsthm}
\usepackage{setspace}
\usepackage{multirow}
\usepackage{hyperref}
\usepackage{cases}
\usepackage[bottom]{footmisc}
\usepackage[small]{caption}
\usepackage[usenames,dvipsnames]{pstricks}
\usepackage{epsfig}
\usepackage{pst-grad} 
\usepackage{pst-plot} 

\setmarginsrb{1.8cm}{1.8cm}{1.8cm}{1.8cm}{0cm}{0cm}{0cm}{0.5cm}

\begin{document}

\newtheorem{Proposition}{Proposition}
\newtheorem{Theorem}{Theorem}
\newtheorem{Corollary}{Corollary}
\newtheorem{Definition}{Definition}
\newtheorem{Lemma}{Lemma}

\title{From infinity to one: The reduction of some mean field games to a global control problem\addtocounter{footnote}{-1}\thanks{The author wishes to acknowledge the helpful conversations with Yves Achdou (Université Paris-Diderot), Guy Barles (Université de Tours), Diogo Gomes (IST Lisboa), Jean-Michel Lasry (Université Paris-Dauphine) and Pierre-Louis Lions (Collège de France). Also, the research presented in this paper benefited from the support of the Conseil Français de l'Energie.}}
\author{Olivier Guéant\addtocounter{footnote}{0}\thanks{UFR de Math\'ematiques, Laboratoire Jacques-Louis Lions, Universit\'e Paris-Diderot. 175, rue du Chevaleret, 75013 Paris, France. \texttt{olivier.gueant@ann.jussieu.fr}}}
\date{Preliminary version, Summer 2011}

\maketitle
\abstract{This paper presents recent results from Mean Field Game theory underlying the introduction of common noise that imposes to incorporate the distribution of the agents as a state variable. Starting from the usual mean field games equations introduced in \cite{MFG1,MFG2,MFG3} and adapting them to games on graphs, we introduce a partial differential equation, often referred to as the Master equation (see \cite{MFG4}), from which the MFG equations can be deduced. Then, this Master equation can be reinterpreted using a global control problem inducing the same behaviors as in the non-cooperative initial mean field game.}

\section*{Introduction}

Mean field games have been introduced in 2006 by J.-M. Lasry and P.-L. Lions \cite{MFG1,MFG2,MFG3} to generalize stochastic differential games to games involving a continuum of players.\\
Since 2006, many applications of mean field games have been proposed, particularly in economics (see for instance \cite{gueant2008mean}, \cite{ParisPrinceton}, or \cite{lachapelle2010computation}) and an important effort has been made to solve the partial differential equations associated to mean field games when both time and the state space are continuous (see \cite{achdou2010capuzzo}, \cite{achdou2010mean}, \cite{gueant2011}, \cite{lachapelle2010computation}, etc).\\
In this paper, we consider a mean field game on a complex discrete state space that is represented by a graph and we explain the discrete counterpart of the MFG equations and of the Master equation formally introduced in \cite{MFG4}. This Master equation, or in our discrete case the Master equations, ``contains'' the MFG equations and allows for the introduction of random noise (although this question is not tackled here). In addition to these Master equations that can be exposed in a simple way in this discrete context without relying on differential calculus in a space of probability measures, we provide an interpretation of them as the equations solved by the partial differentials of a function that solves a single Hamilton-Jacobi equation in high dimension associated to a global optimization problem on the graph.\\

In the first section, we introduce the setting of our mean field game on a graph and define what a Nash-MFG equilibrium is. Then, in section 2, we define the MFG equations in the case of a graph, along with the Master equations, and we show how to deduce a Nash-MFG equilibrium from smooth solutions of these equations. In section 3 we define a global optimization problem on the graph and show how the solution of the associated Hamilton-Jacobi equation is related to the Master equations and eventually to the initial mean field game.

\section{Mean field games in a graph}

We consider a directed graph $\mathcal{G}$ whose nodes are indexed by integers from 1 to $N$. For each node $i \in \mathcal{N}= \lbrace 1, \ldots, N\rbrace$ we introduce $\mathcal{V}(i) \subset \mathcal{N}\setminus \lbrace i\rbrace$ the set of nodes $j$ for which a directed edge exists from $i$ to $j$. The cardinal of this set is denoted $d_i$ and called the out-degree of the node $i$. Reciprocally, we denote $\mathcal{V}^{-1}(i) \subset \mathcal{N}\setminus \lbrace i\rbrace$ the set of nodes $j$ for which a directed edge exists from $j$ to $i$.\\
This graph $\mathcal{G}$ is going to be the state space of our mean field game and we suppose that there is a continuum of anonymous and identical players of size $1$, playing a game on the graph as described in the next paragraphs\footnote{While writing this article, we learnt during a conference organized in Rome on Mean Field Games that results on the limit of continuous time and finite state space $M$-player-games as $M \to +\infty$ had been obtained by Diogo Gomes, justifying our direct passage to a continuum of players}.\\

Each player's position is represented by a continuous-time Markov chain $(X_t)_t$ with values in $\mathcal{G}$ and instantaneous transition probabilities at time $t$ described by a collection of functions $\lambda_t(i,\cdot): \mathcal{V}(i) \to \mathbb{R}_+$ (for each node $i \in \mathcal{N}$).\\
Each player is able to decide on the values of the transition probabilities at time $t$ and pays an instantaneous cost\footnote{There is no additional difficulty when a time dependency is added.} $L(i,(\lambda_{i,j})_{j \in \mathcal{V}(i)})$ to set the value of $\lambda(i,j)$ to $\lambda_{i,j}$.\\

The assumptions made on the continuous functions $L(i,\cdot)$ are the following.\\
First, we assume asymptotic super-linearity for each $i \in \mathcal{N}$:
$$
 \forall i \in \mathcal{N}, \lim_{\lambda \in \mathbb{R}_+^{d_i}, |\lambda| \to +\infty} \frac{L(i,\lambda)}{|\lambda|} = + \infty
$$
Then, we define the Hamiltonians $H(i,\cdot)$:
$$
 \forall i \in \mathcal{N}, p \in \mathbb{R}^{d_i} \mapsto H(i,p) = \sup_{\lambda \in \mathbb{R}_+^{d_i}} \lambda\cdot p - L(i,\lambda)
$$
and we assume that for each $i\in \mathcal{N}$, $H(i,\cdot)$ is a $C^1$ function with $$\nabla H(i,p) = \mathrm{argmax}_{\lambda \in \mathbb{R}_+^{d_i}} \lambda\cdot p - L(i,\lambda)$$
This assumption is satisfied as soon as for each $i\in \mathcal{N}$, $L(i,\cdot)$ is a $C^2$ strongly convex function or if $L(i,\lambda) = \sum_{j\in\mathcal{V}(i)} C(i,j,\lambda_{i,j})$ with each of the functions $C(i,j,\cdot)$ a strictly convex $C^2$ function.\\

Now, to define the mean field game, let us introduce a domain $\Omega \subset \mathbb{R}^N$ containing $\mathcal{P}_N = \lbrace (x_1, \ldots, x_N) \in [0,1]^n , \sum_{i=1}^N x_i = 1\rbrace$ and two functions $f$ and $g$ defined\footnote{If one adds a time dependency to $f$, the following results can be adapted straightforwardly.} on $\mathcal{N}\times \Omega$ such that $\forall i \in \mathcal{N}$, $f(i,\cdot)$ and $g(i,\cdot)$ are continuous on $\Omega$.

Now, let us also define the set $\mathcal{A}$ of admissible controls. We are only considering markovian controls and we define:
$$\mathcal{A} = \left\lbrace (\lambda_t(i,j))_{t \in [0,T],i \in \mathcal{N}, j\in \mathcal{V}(i)} \mathrm{\; deterministic} \left| \forall i \in \mathcal{N}, \forall j \in \mathcal{V}(i), t \mapsto \lambda_t(i,j) \in L^{\infty}(0,T)\right. \right\rbrace$$

Now, for a given admissible control\footnote{We call $\lambda$ a control although it is an abuse of terminology since the controls consist in the values of $\lambda$.} $\lambda \in \mathcal{A}$ and a given function $m: t \in [0,T] \mapsto (m(t,1),\ldots,m(t,N)) \in \mathcal{P}_N$ we define the payoff function $J_m : [0,T]\times\mathcal{N}\times\mathcal{A} \to \mathbb{R}$ by:

$$J_m(t,i,\lambda) = \mathbb{E}\left[\int_{t}^T \left(-L(X_s,\lambda_s(X_s,\cdot)) + f(X_s,m(s))\right) ds + g\left(X_T,m(T)\right)\right] $$

for $(X_s)_{s \in [t,T]}$ a Markov chain on $\mathcal{G}$, starting from $i$ at time $t$, with instantaneous transition probabilities given by $(\lambda_s)_{s \in [t,T]}$.\\

From this we can adapt the definition of a (symmetric) Nash-MFG equilibrium to our context:

\begin{Definition}[Nash-MFG symmetric equilibrium]
A differentiable function $m: t \in [0,T] \mapsto (m(t,1),\ldots,m(t,N)) \in \mathcal{P}_N$ is said to be a Nash-MFG equilibrium, if there exists an admissible control $\lambda \in \mathcal{A}$ such that:
$$\forall \tilde{\lambda} \in \mathcal{A}, \forall i \in \mathcal{N} , J_m(0,i,\lambda) \ge J_m(0,i,\tilde{\lambda})$$
and
$$\forall i \in \mathcal{N}, \frac{d\ }{dt} m(t,i) = \sum_{j \in \mathcal{V}^{-1}(i)} \lambda_t(j,i) m(t,j) - \sum_{j \in \mathcal{V}(i)} \lambda_t(i,j) m(t,i)$$
In that case, $\lambda$ is called an optimal control.\\
\end{Definition}

In this definition, $m$ stands implicitly for the distribution of the players across the different nodes.\\

As for all mean field games, the problem is a fixed point problem since, starting from a (guessed) trajectory $m$, the control problem in Definition 1 induces trajectories for the players and to obtain a Nash-MFG equilibrium these individual trajectories must generate the right macroscopic trajectory $m$ for the distribution of the players in $\mathcal{G}$.\\

In the next section we are going to write the MFG equations associated to our problem, which we call the \emph{decentralized problem} as it involves all the players, and see that these equations can be extracted from a set of equations called the Master equations, as introduced by J.-M. Lasry and P.-L. Lions and presented in \cite{MFG4} in a different framework.

\section{MFG equations and the Master equation}

As for any control problem we can associate, when $m$ is fixed, a Hamilton-Jacobi equation to the above control problem. More exactly, since there are $N$ states, we need to write down $N$ Hamilton-Jacobi equations, one for each state. However, these backward equations need to be coupled with a system of forward transport equations that describe the resulting evolution of the distribution of players.\\

The graph counterpart of the two partial differential equations of classical mean field games (see \cite{MFG1,MFG2,MFG3}) are the following:

\begin{Definition}[The $\mathcal{G}$-MFG equations]
The $\mathcal{G}$-MFG equations consist in a system of $2N$ equations, the unknown being $ t \in [0,T] \mapsto (u(t,1), \ldots, u(t,N), m(t,1), \ldots, m(t,N))$.\\

The first half of these equations are the Hamilton-Jacobi equations associated to the above problem and consist in the following system:
$$\forall i \in \mathcal{N}, \quad \frac{d\ }{dt} u(t,i) + H\left(i,(u(t,j)-u(t,i))_{j \in \mathcal{V}(i)}\right) + f(i,m(t,1),\ldots,m(t,N)) =0$$
with $u(T,i) = g(i,m(t,1),\ldots,m(t,N))$.\\

The second half of these equations are forward transport equations:
$$\forall i \in \mathcal{N}, \quad \frac{d\ }{dt} m(t,i) = \sum_{j \in \mathcal{V}^{-1}(i)} \frac{\partial H(j,\cdot)}{\partial{p_i}}\left((u(t,k)-u(t,j))_{k \in \mathcal{V}(j)}\right) m(t,j)$$$$ - \sum_{j \in \mathcal{V}(i)} \frac{\partial H(i,\cdot)}{\partial{p_j}}\left((u(t,k)-u(t,i))_{k \in \mathcal{V}(i)}\right) m(t,i)$$
with $(m(0,1), \ldots, m(0,N)) = m^0 \in \mathcal{P}_N$ given.\\
\end{Definition}

Now, we enounce a result saying that a $C^1$ solution of the $\mathcal{G}$-MFG equations provides a Nash-MFG equilibrium and an optimal control:

\begin{Proposition}[The $\mathcal{G}$-MFG equations as a sufficient condition]
Let $m^0 \in \mathcal{P}_N$ and let us consider a $C^1$ solution $(u(t,1), \ldots, u(t,N), m(t,1), \ldots, m(t,N))$ of the $\mathcal{G}$-MFG equations with $(m(0,1), \ldots, m(0,N)) = m^0$.\\

Then $t \mapsto m(t) = (m(t,1), \ldots, m(t,N))$ is a Nash-MFG equilibrium and the relations $\lambda_t(i,j) = \frac{\partial H(i,\cdot)}{\partial{p_j}}\left((u(t,k)-u(t,i))_{k \in \mathcal{V}(i)}\right)$ define an optimal control.
\end{Proposition}

\textit{Proof:}\\

Let us consider an admissible control $\tilde{\lambda} \in \mathcal{A}$. Let us then consider an initial position $i\in \mathcal{N}$ and a Markov chain $(X^{0,i,\tilde{\lambda}}_t)_t$ starting in $i$ at time $0$, with instantaneous transition probabilities given by $\tilde{\lambda}$.\\
Since the values taken by $\lambda$ are bounded, there is almost surely a finite number $Q$ of moves in the time interval $[0,T)$ and we can write (almost surely):

$$u\left(T,X^{0,i,\tilde{\lambda}}_{T-}\right) = u(0,i) + \left[\sum_{q=1}^Q u\left(t_q,X^{0,i,\tilde{\lambda}}_{t_q}\right) - u\left(t_{q-1},X^{0,i,\tilde{\lambda}}_{t_{q-1}}\right)\right] + u\left(T,X^{0,i,\tilde{\lambda}}_{T-}\right) - u\left(t_Q,X^{0,i,\tilde{\lambda}}_{t_Q}\right)$$
where $t_0 = 0$ and the other $t_i$s are the time of jumps.\\

Then:
$$u\left(T,X^{0,i,\tilde{\lambda}}_{T-}\right) = u(0,i) + \int_0^T \frac{\partial u}{\partial t}(t,X_{t-}) dt + \sum_{q=1}^Q u\left(t_q,X^{0,i,\tilde{\lambda}}_{t_q}\right) - u\left(t_{q},X^{0,i,\tilde{\lambda}}_{t_{q}-}\right) $$
$$=u(0,i) + \int_0^T \frac{\partial u}{\partial t}(t,X_{t-}) dt + \int_0^T \sum_{j \in \mathcal{V}(X_{t-})} \left(u\left(t,j\right) - u\left(t,X^{0,i,\tilde{\lambda}}_{t-}\right)\right) \tilde{\lambda}_t(X_{t-},j) dt + M_T$$
where $(M_t)_t$ is a local martingale starting from $0$ that happens to be a martingale since $u$ is bounded.\\

Hence:

$$\mathbb{E}\left[g(X^{0,i,\tilde{\lambda}}_T,m(T))\right] = \mathbb{E}\left[g(X^{0,i,\tilde{\lambda}}_{T-},m(T-))\right] = u(0,i)$$$$ + \mathbb{E}\left[\int_0^T \left(\frac{\partial u}{\partial t}(t,X_{t-}) + \sum_{j \in \mathcal{V}(X_{t-})} \left(u\left(t,j\right) - u\left(t,X^{0,i,\tilde{\lambda}}_{t-}\right)\right) \tilde{\lambda}_t(X_{t-},j)\right) dt \right]$$

$$\implies\mathbb{E}\left[g(X^{0,i,\tilde{\lambda}}_T,m(T)) - \int_{0}^T L(X_{t-},\lambda_t(X_{t-},\cdot)) dt \right]= u(0,i) $$$$+ \mathbb{E}\left[\int_0^T \left(\frac{\partial u}{\partial t}(t,X_{t-}) -L(X_{t-},\lambda_t(X_{t-},\cdot)) + \sum_{j \in \mathcal{V}(X_{t-})} \left(u\left(t,j\right) - u\left(t,X^{0,i,\tilde{\lambda}}_{t-}\right)\right) \tilde{\lambda}_t(X_{t-},j)\right) dt \right]$$

By definition of $u$, we have that:

$$\mathbb{E}\left[g(X^{0,i,\tilde{\lambda}}_T,m(T)) + \int_{0}^T \left(-L(X_{t-},\lambda_t(X_{t-},\cdot))\right) dt \right]\le u(0,i) - \mathbb{E}\left[\int_{0}^T f(X^{0,i,\tilde{\lambda}}_{t-},m(t-)) dt \right]$$
$$\implies \mathbb{E}\left[\int_{0}^T \left(-L(X_{t},\lambda_t(X_{t},\cdot))+ f(X^{0,i,\tilde{\lambda}}_{t},m(t))\right) dt + g(X^{0,i,\tilde{\lambda}}_T,m(T)) \right]\le u(0,i)$$

with equality when $\tilde{\lambda}$ is taken to be the control $\lambda$ as described in the statement of the result (because of our assumption on the hamiltonians), this control being indeed admissible since the function $u$ is bounded and since each $H(i, \cdot)$ is $C^1$.\\

Hence:
$$J_m(0,i,\tilde{\lambda}) \le u(0,i) = J_m(0,i,\lambda)$$
and this proves the result.\qed\\

Existence and uniqueness of solutions to the $\mathcal{G}$-MFG equations are now tackled. We first start with existence and our proof is based on a Schauder fixed-point argument and a priori estimates to obtain compactness. Then we will present a criterion to ensure uniqueness of $C^1$ solutions.\\

For the existence result, we first start with a Lemma stating that, for a fixed $m$, the $N$ Hamilton-Jacobi equations amongst the $\mathcal{G}$-MFG equations obey a comparison principle:

\begin{Lemma}[Comparison principle]
Let $m : [0,T] \to \mathcal{P}_N$ be a continuous function.
Let $u : t \in [0,T] \mapsto (u(t,1), \ldots, u(t,N))$ be a $C^1$ function that verifies:
$$\forall i \in \mathcal{N}, \quad -\frac{d\ }{dt} u(t,i) - H\left(i,(u(t,j)-u(t,i))_{j \in \mathcal{V}(i)}\right) - f(i,m(t,1),\ldots,m(t,N)) \le0$$
with $u(T,i) \le g(i,m(T,1),\ldots,m(T,N))$.\\

Let $v : t \in [0,T] \mapsto (v(t,1), \ldots, v(t,N))$ be a $C^1$ function that verifies:
$$\forall i \in \mathcal{N}, \quad -\frac{d\ }{dt} v(t,i) - H\left(i,(v(t,j)-v(t,i))_{j \in \mathcal{V}(i)}\right) - f(i,m(t,1),\ldots,m(t,N)) \ge0$$
with $v(T,i) \ge g(i,m(T,1),\ldots,m(T,N))$.\\

Then, $\forall i \in \mathcal{N}, \forall t \in [0,T], v(t,i) \ge u(t,i)$.\\
\end{Lemma}

\textit{Proof:}\\

Let us consider for a given $\epsilon > 0$, a point $(t^*,i^*) \in [0,T]\times\mathcal{N}$ such that $$u(t^*,i^*) - v(t^*,i^*) - \epsilon (T-t^*) = \max_{(t,i) \in [0,T]\times\mathcal{N}}  u(t,i) - v(t,i) - \epsilon (T-t)$$

If $t^*\in [0,T)$, then $\left.\frac{d\ }{dt} \left(u(t,i^*) - v(t,i^*) - \epsilon (T-t)\right)\right|_{t=t^*} \le 0$. Also, by definition of $(t^*,i^*)$, $\forall j \in \mathcal{V}(i^*)$, $u(t^*,i^*) - v(t^*,i^*) \ge u(t^*,j) - v(t^*,j)$ and hence, by definition of $H(i^*,\cdot)$:
$$H\left(i^*, (v(t^*,j)-v(t^*,i^*))_{j \in \mathcal{V}(i^*)}\right) \ge H\left(i^*, (u(t^*,j)-u(t^*,i^*))_{j \in \mathcal{V}(i^*)}\right) $$

Combining these inequalities we get:
$$ -\frac{d\ }{dt} u(t^*,i^*) - H\left(i^*,(u(t^*,j)-u(t^*,i^*))_{j \in \mathcal{V}(i^*)}\right) \ge -\frac{d\ }{dt} v(t^*,i^*) - H\left(i^*,(v(t^*,j)-v(t^*,i^*))_{j \in \mathcal{V}(i^*)}\right) + \epsilon$$

But this is in contradiction with the hypotheses on $u$ and $v$.\\

Hence $t^*=T$ and $\max_{(t,i) \in [0,T]\times\mathcal{N}}  u(t,i) - v(t,i) - \epsilon (T-t) \le 0$ because of the assumptions on $u(T,i)$ and $v(T,i)$.\\

This being true for any $\epsilon > 0$, we have that $\max_{(t,i) \in [0,T]\times\mathcal{N}}  u(t,i) - v(t,i) \le 0$.\qed\\

This lemma allows to provide a bound to any solution $u$ of the $N$ Hamilton-Jacobi equations and this bound is then used to obtain compactness in order to apply Schauder's fixed point theorem.

\begin{Proposition}[Existence of a solution to the $\mathcal{G}$-MFG equations]
Let $m^0 \in \mathcal{P}_N$. Under the assumptions made in section 1, there exists a $C^1$ solution $(u,m)$ of the $\mathcal{G}$-MFG equations such that $m(0) = m^0$.
\end{Proposition}

\textit{Proof:}\\

Let $m : [0,T] \to \mathcal{P}_N$ be a continuous function.\\

Let then consider the solution $u : t \in [0,T] \mapsto (u(t,1), \ldots, u(t,N))$ to the Hamilton-Jacobi equations:

$$\forall i \in \mathcal{N}, \quad \frac{d\ }{dt} u(t,i) + H\left(i,(u(t,j)-u(t,i))_{j \in \mathcal{V}(i)}\right) + f(i,m(t,1),\ldots,m(t,N)) = 0$$
with $u(T,i) = g(i,m(T,1),\ldots,m(T,N))$.\\

This function $u$ is a well defined $C^1$ function with the following bound coming from the above lemma:

$$\sup_{i \in \mathcal{N}} \|u(\cdot,i)\|_{\infty} \le \sup_{i \in \mathcal{N}} \|g(i,\cdot)\|_{\infty} + \left(\sup_{i \in \mathcal{N}} \|f(i,\cdot)\|_{\infty} + \sup_{i \in \mathcal{N}} |H(i,0)|\right)T$$

Using this bound and the fact that the functions $H(i,\cdot)$ are assumed to be $C^1$, we can define a function $\tilde{m} : [0,T] \to \mathcal{P}_N$ by:

$$\forall i \in \mathcal{N}, \quad \frac{d\ }{dt} \tilde{m}(t,i) = \sum_{j \in \mathcal{V}^{-1}(i)} \frac{\partial H(j,\cdot)}{\partial{p_i}}\left((u(t,k)-u(t,j))_{k \in \mathcal{V}(j)}\right) \tilde{m}(t,j)$$$$ - \sum_{j \in \mathcal{V}(i)} \frac{\partial H(i,\cdot)}{\partial{p_j}}\left((u(t,k)-u(t,i))_{k \in \mathcal{V}(i)}\right) \tilde{m}(t,i)$$
with $(\tilde{m}(0,1), \ldots, \tilde{m}(0,N)) = m^0 \in \mathcal{P}_N$.\\

$\frac{d\tilde{m}}{dt}$ is bounded, the bounds depending only on the functions $f(i,\cdot)$, $g(i,\cdot)$ and $H(i,\cdot)$, $i \in \mathcal{N}$.\\

As a consequence, if we define $\Theta : m \in C([0,T],\mathcal{P}_N) \mapsto \tilde{m} \in C([0,T],\mathcal{P}_N)$, $\Theta$ is a continuous function (from classical ODEs theory) with $\Theta(C([0,T],\mathcal{P}_N))$ a relatively compact set (because of Ascoli's Theorem and the uniform Lipschitz property we just obtained).\\

Hence, because $C([0,T],\mathcal{P}_N)$ is convex, by Schauder's fixed point theorem, there exists a fixed point $m$ to $\Theta$. If we then consider $u$ associated to $m$ by the Hamilton-Jacobi equations as above, $(u,m)$ is a $C^1$ solution to the $\mathcal{G}$-MFG equations.\qed\\

Coming now to uniqueness, we deduce from the proof of uniqueness introduced in \cite{MFG1,MFG2,MFG3}, a criterion to ensure uniqueness of $C^1$ solutions of the $\mathcal{G}$-MFG equations:

\begin{Proposition}[Uniqueness for the solution of the $\mathcal{G}$-MFG equations]
Assume that $f$ and $g$ are such that:
$$\forall (m^1,m^2) \in \mathcal{P}_N \times \mathcal{P}_N, \sum_{i=1}^N (f(i,m^1_1,\ldots,m^1_N) - f(i,m^2_1,\ldots,m^2_N))(m^1_i - m^2_i) \ge 0 \implies m^1=m^2$$
and
$$\forall (m^1,m^2) \in \mathcal{P}_N \times \mathcal{P}_N, \sum_{i=1}^N (g(i,m^1_1,\ldots,m^1_N) - g(i,m^2_1,\ldots,m^2_N))(m^1_i - m^2_i) \ge 0 \implies m^1=m^2$$
Then, if $(u,m)$ and $(\tilde{u},\tilde{m})$ are two $C^1$ solutions of the $\mathcal{G}$-MFG equations, we have $m=\tilde{m}$ and $u=\tilde{u}$.
\end{Proposition}

\textit{Proof:}\\

The proof of this result is a straightforward adaptation of the classical proof of uniqueness for continuous MFG equations.\\

It consists in computing the value of $I = \int_0^T \sum_{i=1}^N \frac{d\ }{dt} \left((u(t,i) - \tilde{u}(t,i))(m(t,i) - \tilde{m}(t,i))\right) dt $ in two different ways.\\

We first know directly that $$I = \sum_{i=1}^N (g(i,m(T,1),\ldots,m(T,N)) - g(i,\tilde{m}(T,1),\ldots,\tilde{m}(T,N)))(m(T,i) - \tilde{m}(T,i))$$

Now, differentiating the product we get, after reordering the terms that:

$$I = -\int_0^T \sum_{i=1}^N (f(i,m(t,1),\ldots,m(t,N)) - f(i,\tilde{m}(t,1),\ldots,\tilde{m}(t,N)))(m(t,i) - \tilde{m}(t,i)) dt$$
$$+ \int_0^T  \sum_{i=1}^N m(t,i) \Bigg{\lbrack} H(i,(\tilde{u}(t,k)-\tilde{u}(t,i))_{k \in \mathcal{V}(i)}) - H(i,(u(t,k)-u(t,i))_{k \in \mathcal{V}(i)})$$$$ - \sum_{j \in \mathcal{V}(i)} ((u(t,i) - \tilde{u}(t,i)) - (u(t,j) - \tilde{u}(t,j)) ) \frac{\partial H(i,\cdot)}{\partial{p_j}}\left((u(t,k)-u(t,i))_{k \in \mathcal{V}(i)}\right) \Bigg{\rbrack}dt$$
$$+\int_0^T  \sum_{i=1}^N \tilde{m}(t,i) \Bigg{\lbrack} H(i,(u(t,k)-u(t,i))_{k \in \mathcal{V}(i)}) - H(i,(\tilde{u}(t,k)-\tilde{u}(t,i))_{k \in \mathcal{V}(i)}) $$$$ - \sum_{j \in \mathcal{V}(i)} ((\tilde{u}(t,i) - u(t,i)) - (\tilde{u}(t,j) - u(t,j)) ) \frac{\partial H(i,\cdot)}{\partial{p_j}}\left((\tilde{u}(t,k)-\tilde{u}(t,i))_{k \in \mathcal{V}(i)}\right) \Bigg{\rbrack}dt$$

Now, using the convexity of the hamiltonian functions, we know that:

$$\int_0^T \sum_{i=1}^N (f(i,m(t,1),\ldots,m(t,N)) - f(i,\tilde{m}(t,1),\ldots,\tilde{m}(t,N)))(m(t,i) - \tilde{m}(t,i)) dt$$
$$+\sum_{i=1}^N (g(i,m(T,1),\ldots,m(T,N)) - g(i,\tilde{m}(T,1),\ldots,\tilde{m}(T,N)))(m(T,i) - \tilde{m}(T,i)) \ge 0$$

Hence, using the hypotheses on $f$ and $g$, $m=\tilde{m}$.\\
Using the maximum principle stated in Lemma 1, we hence have $u=\tilde{u}$ and the result is proved.\qed\\

A simple example of functions $f$ and $g$ satisfying the above criterion consists in taking for $f(i,\cdot)$ and $g(i,\cdot)$ functions of $m_i$ only, strictly decreasing in $m_i$.\\

Now we define the graph counterpart of the Master equation proposed by J.-M. Lasry and P.-L. Lions and presented in \cite{MFG4}:

\begin{Definition}[The $\mathcal{G}$-Master equations]
The $\mathcal{G}$-Master equations consist in $N$ equations, the unknown being $ (t,m_1,\ldots,m_N) \in [0,T]\times\Omega \mapsto (U_1(t,m_1,\ldots,m_N), \ldots, U_N(t,m_1,\ldots,m_N))$.\\
$$\forall i \in \mathcal{N}, \quad\frac{\partial{U_i}}{\partial t}(t,m_1, \ldots, m_N) + H\left(i,(U_j(t,m_1, \ldots, m_N)-U_i(t,m_1, \ldots, m_N))_{j \in \mathcal{V}(i)}\right) $$
$$+ \sum_{l=1}^N \frac{\partial{U_i}}{\partial{m_l}}(t,m_1, \ldots, m_N) \left[\sum_{j \in \mathcal{V}^{-1}(l)} \frac{\partial H(j,\cdot)}{\partial{p_l}}\left((U_k(t,m_1, \ldots, m_N)-U_j(t,m_1, \ldots, m_N))_{k \in \mathcal{V}(j)}\right) m_j\right.$$$$\left. - \sum_{j \in \mathcal{V}(l)} \frac{\partial H(l,\cdot)}{\partial{p_j}}\left((U_k(t,m_1, \ldots, m_N)-U_l(t,m_1, \ldots, m_N))_{k \in \mathcal{V}(l)}\right) m_l\right]+  f(i, m_1, \ldots, m_N)  =0$$
with $U_i(T,m_1, \ldots, m_N) = g(i,m_1,\ldots,m_N)$.\\
\end{Definition}

This equation is the Hamilton-Jacobi equation associated to the decentralized problem when we consider that the state variable is not only the current position but also the repartition of the players on $\mathcal{G}$. In particular, when there is a common noise affecting for instance the functions $f$ and $g$, one cannot write the $\mathcal{G}$-MFG equations and one has to rely instead on the $\mathcal{G}$-Master equations (see \cite{MFG4}).\\

To prove that these $\mathcal{G}$-Master equations are indeed more general than the $\mathcal{G}$-MFG equations we enounce the following proposition:

\begin{Proposition}[From $\mathcal{G}$-Master equations to $\mathcal{G}$-MFG equations]
Let us consider a $C^1$ function $(t,m_1, \ldots, m_N) \in [0,T]\times\Omega \mapsto (U_1(t,m_1,\ldots,m_N),\ldots,U_N(t,m_1,\ldots,m_N))$ verifying the $\mathcal{G}$-Master equations.
Let us then consider a $C^1$ function $m: t \in [0,T] \mapsto (m(t,1),\ldots,m(t,N))$ verifying:
$$\forall i \in \mathcal{N}, \quad\frac{d }{dt} m(t,i) =$$$$ \sum_{j \in \mathcal{V}^{-1}(i)} \frac{\partial H(j,\cdot)}{\partial{p_i}}\left((U_k(t,m(t,1), \ldots, m(t,N))-U_j(t,m(t,1), \ldots, m(t,N)))_{k \in \mathcal{V}(j)}\right) m(t,j)$$$$ - \sum_{j \in \mathcal{V}(i)} \frac{\partial H(i,\cdot)}{\partial{p_j}}\left((U_k(t,m(t,1), \ldots, m(t,N))-U_i(t,m(t,1), \ldots, m(t,N)))_{k \in \mathcal{V}(i)}\right) m(t,i)$$
with $(m(0,1), \ldots, m(0,N)) = m^0 \in \mathcal{P}_N$.\\
If we define $$\forall i \in \mathcal{N}, u(t,i) = U_i(t,m(t,1), \ldots, m(t,N))$$
then $t \in [0,T] \mapsto (u(t,1), \ldots, u(t,N), m(t,1), \ldots, m(t,N))$ is a solution of the $\mathcal{G}$-MFG equations with the initial condition $m^0$.
\end{Proposition}

\textit{Proof:}\\

Differentiating $u(t,i) = U_i(t,m(t,1), \ldots, m(t,N))$ with respect to $t$, one obtain the $i^{th}$ Hamilton-Jacobi equation that defines the $\mathcal{G}$-MFG equations. The result is then straightforward.\qed\\

\section{Potential games and reduction to a planning problem}

In this section, we provide a way to deduce the $\mathcal{G}$-Master equations that consist in $N$ equations from a single equation when the functions $f$ and $g$ verify the hypotheses of potential games\footnote{This is for instance the case if $\forall i\in \mathcal{N}$, $f(i,\cdot)$ and $g(i,\cdot)$ only depend on $m_i$.}. We say that the above game is a potential game when there exist two $C^1$ functions $F: (m_1,\ldots, m_N) \in \Omega \mapsto F(m_1,\ldots, m_N)$ and $G: (m_1,\ldots, m_N) \in \Omega \mapsto G(m_1,\ldots, m_N)$ such that $\forall i \in \mathcal{N}$, $\frac{\partial F}{\partial{m_i}} = f(i,\cdot)$ and $ \frac{\partial G}{\partial{m_i}} = g(i,\cdot)$ in $\Omega$.\\

To that purpose let us introduce a global optimization problem that we phrase a \emph{planning problem}\footnote{The term ``\emph{planning problem}'' must be understood here as the problem of a planner controlling the whole population of players. The planing problem studied in \cite{achdou2010capuzzo} is of a different nature.}. This optimization problem consists in optimizing directly, with a macroscopic viewpoint, the different flows of agents (the preceding players being now controlled by a single agent, often referred to as a planner) between the nodes of the graph.\\

We introduce for $t \in [0,T]$, $m^t \in \mathcal{P}_N$ and a given admissible control (function) $\lambda \in \mathcal{A}$, the payoff function

$$\mathcal{J}(t,m^t,\lambda) = \int_t^T \left(F(m(s)) - \sum_{i=1}^N L(i,(\lambda_s(i,j))_{j\in \mathcal{V}(i)}) m(s,i) \right) ds + G(m(T))$$
where $\forall i \in \mathcal{N}, m(t,i) = m^t_i$ and
$$\forall i \in \mathcal{N}, \forall s \in [t,T], \frac{d\ }{ds} m(s,i) = \sum_{j \in \mathcal{V}^{-1}(i)} \lambda_s(j,i) m(s,j) - \sum_{j \in \mathcal{V}(i)} \lambda_s(i,j) m(s,i)$$

The optimization problem we consider is, for a given $m^0 \in \mathcal{P}_N$:

$$\sup_{\lambda \in \mathcal{A}} \mathcal{J}(0,m^0,\lambda)$$

Noteworthy, and contrary to the first optimization problem that was a stochastic control problem, this new problem is a deterministic control problem.\\

We define the Hamilton-Jacobi equation associated to this deterministic control problem:

\begin{Definition}[The $\mathcal{G}$-planning equation]
The $\mathcal{G}$-planning equation consists in a single partial differential equation, the unknown being $ (t,m) \in [0,T]\times\Omega \mapsto \Phi(t,m_1,\ldots,m_N)$:\\

$$\frac{\partial \Phi}{\partial t}(t,m_1, \ldots,m_N) + \mathcal{H}(m_1,\ldots,m_N,\nabla \Phi) + F(m_1, \ldots,m_N) = 0 $$
with the terminal conditions $\Phi(T,m_1,\ldots,m_N) = G(m_1,\ldots,m_N)$, where the hamiltonian is given by:

$$\mathcal{H}(m_1,\ldots,m_N,p) = \sup_{(\lambda_{i,j})_{i \in \mathcal{N}, j\in \mathcal{V}(i)}}\sum_{i=1}^N\left[\left(\sum_{j \in \mathcal{V}^{-1}(i)} \lambda_{j,i} m_j - \sum_{j \in \mathcal{V}(i)} \lambda_{i,j} m_i\right)p_i  - L(i,(\lambda_{i,j})_{j\in \mathcal{V}(i)})m_i\right]$$
\end{Definition}

Our goal is to link this \emph{planning problem} to the \emph{decentralized problem}. The first step is to write the hamiltonian $\mathcal{H}$ in another way:

\begin{Lemma}
$$\mathcal{H}(m_1,\ldots,m_N,p) = \sum_{i=1}^N m_i H\left(i,\left(p_j - p_i\right)_{j \in \mathcal{V}(i)}\right)$$
\end{Lemma}

\emph{Proof:}\\

Reordering the terms we can write $\mathcal{H}(m_1,\ldots,m_N,p)$ as:
$$\sum_{i=1}^N m_i \sup_{(\lambda_{i,j})_{j\in \mathcal{V}(i)}} \left(\sum_{j \in \mathcal{V}(i)} \lambda_{i,j} (p_j-p_i)  - L(i,(\lambda_{i,j})_{j\in \mathcal{V}(i)})\right)$$
Then using the definition of the hamiltonians in the decentralized problem, we get the result.\qed\\

\begin{Proposition}
Let us consider a $C^1$ function $\Phi$ solution of the $\mathcal{G}$-planning equation. Then, $\Phi$ restricted to $[0,T]\times\mathcal{P}_N$ is the value function of the above planning problem, i.e.:
$$\forall (t,m^t) \in [0,T]\times \mathcal{P}_N,\Phi(t,m^t) = \sup_{\lambda \in \mathcal{A}} \mathcal{J}(t,m^t,\lambda)$$
Moreover, if we define $\forall i \in \mathcal{N}, m(t,i) = m^t_i$ and
$$\forall i \in \mathcal{N}, \forall s \in [t,T], \frac{d\ }{ds} m(s,i) = \sum_{j \in \mathcal{V}^{-1}(i)} \lambda_s(j,i) m(s,j) - \sum_{j \in \mathcal{V}(i)} \lambda_s(i,j) m(s,i)$$
with $$\lambda_s(i,j) = \frac{\partial H(i,\cdot)}{\partial p_j}\left(\left(\frac{\partial \Phi}{\partial m_k}(s,m(s,1), \ldots,m(s,N))-\frac{\partial \Phi}{\partial m_i}(s,m(s,1), \ldots,m(s,N))\right)_{k \in \mathcal{V}(i)}\right)$$
then $\lambda$ is an optimal control for the planning problem.
\end{Proposition}

\textit{Proof:}\\

Let us consider an admissible control $\tilde{\lambda} \in \mathcal{A}$ and a repartition $m^t \in \mathcal{P}_N$.\\
Then, let us define:
$$\forall i \in \mathcal{N}, \forall s \in [t,T], \frac{d\ }{ds} m^{t,m^t,\tilde{\lambda}}(s,i) = \sum_{j \in \mathcal{V}^{-1}(i)} \tilde{\lambda}_s(j,i) m^{t,m^t,\tilde{\lambda}}(s,j) - \sum_{j \in \mathcal{V}(i)} \tilde{\lambda}_s(i,j) m^{t,m^t,\tilde{\lambda}}(s,i)$$
with $m(t,i) = m^t_i$.\\
We can write
$$\Phi\left(T,m^{t,m^t,\tilde{\lambda}}(T)\right) - \Phi\left(t,m^t\right) = \int_t^T \frac{\partial \Phi}{\partial s}(s,m^{t,m^t,\tilde{\lambda}}(s)) ds$$
$$+ \sum_{i=1}^N \frac{\partial \Phi}{\partial m_i}(s,m^{t,m^t,\tilde{\lambda}}(s)) \left[\sum_{j \in \mathcal{V}^{-1}(i)} \tilde{\lambda}_s(j,i) m^{t,m^t,\tilde{\lambda}}(s,j) - \sum_{j \in \mathcal{V}(i)} \tilde{\lambda}_s(i,j) m^{t,m^t,\tilde{\lambda}}(s,i)\right]$$
Hence, using the definition of $\Phi$:

$$G\left(m^{t,m^t,\tilde{\lambda}}(T)\right) - \int_t^T\sum_{i=1}^N L(i,(\tilde{\lambda}_s(i,j))_{j\in \mathcal{V}(i)}) m^{t,m^t,\tilde{\lambda}}(s,i) ds$$
$$\le \Phi\left(t,m^t\right) - \int_t^T F(m^{t,m^t,\tilde{\lambda}}(s))ds$$
with equality when $\tilde{\lambda}$ is equal to $\lambda$ as defined in the Proposition, thanks to the preceding Lemma (this $\lambda$ being admissible because of our hypotheses on the hamiltonians, because $\Phi$ is $C^1$ and because $\mathcal{P}_N$ is compact).\\

As a consequence:
$$\mathcal{J}(t,m^t,\tilde{\lambda}) \le \Phi(t,m^t) = \mathcal{J}(t,m^t,\lambda)$$
and this proves the result.\qed\\

Now, we are going to relate the function $\Phi$ to the $\mathcal{G}$-Master equations and prove that a solution to the planning problem can provide a solution to the decentralized problem, that is the initial mean field game.

\begin{Proposition}
Let us consider a $C^2$ function\footnote{In reality, we do not need to be able to differentiate twice with respect to $t$} $\Phi$ solution of the $\mathcal{G}$-planning equation.\\
Define $\forall i \in \mathcal{N}$, $\forall t \in [0,T]$, $\forall m \in \mathcal{P}_N$, $U_i(t,m_1, \ldots, m_N) = \frac{\partial \Phi}{\partial m_i}(t,m_1, \ldots, m_N)$. Then, $\nabla\Phi = U = (U_1, \ldots, U_N)$ verifies the $\mathcal{G}$-Master equations.\\
Consequently, if we define $\forall i \in \mathcal{N}, m(0,i) = m^0_i$ for a given $m^0 \in \mathcal{P}_N$ and
$$\forall i \in \mathcal{N}, \forall s \in [t,T], \frac{d\ }{ds} m(s,i) = \sum_{j \in \mathcal{V}^{-1}(i)} \lambda_s(j,i) m(s,j) - \sum_{j \in \mathcal{V}(i)} \lambda_s(i,j) m(s,i)$$
with $$\lambda_s(i,j) = \frac{\partial H(i,\cdot)}{\partial p_j}\left(\left(\frac{\partial \Phi}{\partial m_k}(s,m(s,1), \ldots,m(s,N))-\frac{\partial \Phi}{\partial m_i}(s,m(s,1), \ldots,m(s,N))\right)_{k \in \mathcal{V}(i)}\right)$$
then $m$ is a Nash-MFG equilibrium and $\lambda$ is an optimal control for the initial mean field game (the decentralized problem).
\end{Proposition}

\textit{Proof:}\\

We need to prove the first assertion of the Proposition, then the other assertions are consequences of Proposition 3.\\

Differentiating the $\mathcal{G}$-planning equation with respect to $m_i$, and using the above Lemma, we know since $\Phi$ is $C^2$ that:

$$\frac{\partial U_i}{\partial t} + H\left(i,(U_j-U_i)_{j \in \mathcal{V}(i)}\right) + \sum_{j=1}^N m_j \sum_{k\in \mathcal{V}(j)} \left(\frac{\partial U_k}{\partial m_i} - \frac{\partial U_j}{\partial m_i} \right) \frac{\partial H(j,\cdot)}{\partial p_k}\left((U_l-U_j)_{l \in \mathcal{V}(j)}\right) + f(i,\cdot) = 0 $$

Now, using the fact that $\Phi$ is $C^2$, we have $\forall (i,j) \in \mathcal{N}^2$, $\frac{\partial U_i}{\partial m_j} = \frac{\partial U_j}{\partial m_i}$. Hence:

$$\frac{\partial U_i}{\partial t} + H\left(i,(U_j-U_i)_{j \in \mathcal{V}(i)}\right) + \sum_{j=1}^N m_j \sum_{k\in \mathcal{V}(j)} \left(\frac{\partial U_i}{\partial m_k} - \frac{\partial U_i}{\partial m_j} \right) \frac{\partial H(j,\cdot)}{\partial p_k}\left((U_l-U_j)_{l \in \mathcal{V}(j)}\right) + f(i,\cdot) = 0 $$

and this is exactly the $i^{th}$ equation amongst the $\mathcal{G}$-Master equations, after reordering the terms. Hence, since the two terminal conditions are coherent, $U = \nabla \Phi$ is indeed solution of the $\mathcal{G}$-Master equations.\qed\\

\section*{Conclusion}

In this paper, we presented a way to deduce a solution to the mean field games equations in a graph from a smooth solution of a single equation Hamilton-Jacobi equation associated to a global deterministic control problem on the whole graph. Moreover, and as described on Figure $1$ below, the Master equations -- that ``contain'' MFG equations and have to be used in the case of games with common noise -- can be reinterpreted as the equations verified by the partial differentials of solutions of this single Hamilton-Jacobi equation.

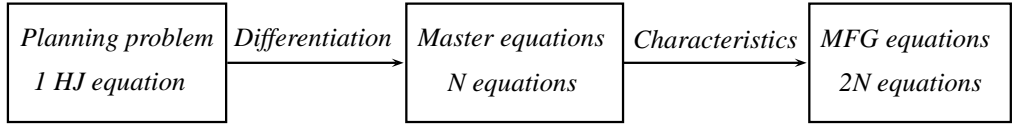
\begin{center}
\scalebox{1} 
{
\begin{pspicture}(0,-0.81)(13.22,0.81)
\psframe[linewidth=0.03,dimen=outer](2.9,0.8)(0.0,-0.8)
\usefont{T1}{ptm}{m}{n}
\rput(1.4510938,0.305){\footnotesize \emph{Planning problem}}
\usefont{T1}{ptm}{m}{n}
\rput(1.350625,-0.275){\footnotesize \emph{1 HJ equation}}
\psline[linewidth=0.03cm,arrowsize=0.05291667cm 2.0,arrowlength=1.4,arrowinset=0.4]{->}(2.9,0)(5.22,0)
\usefont{T1}{ptm}{m}{n}
\rput(4.0265627,0.305){\footnotesize \emph{Differentiation}}
\psframe[linewidth=0.03,dimen=outer](8.12,0.8)(5.24,-0.8)
\usefont{T1}{ptm}{m}{n}
\rput(6.6740627,0.305){\footnotesize \emph{Master equations} }
\usefont{T1}{ptm}{m}{n}
\rput(6.648906,-0.295){\footnotesize \emph{N equations}}
\psline[linewidth=0.03cm,arrowsize=0.05291667cm 2.0,arrowlength=1.4,arrowinset=0.4]{->}(8.1,0)(10.52,0)
\usefont{T1}{ptm}{m}{n}
\rput(9.320937,0.305){\footnotesize \emph{Characteristics}}
\psframe[linewidth=0.03,dimen=outer](13.22,0.8)(10.54,-0.8)
\usefont{T1}{ptm}{m}{n}
\rput(11.864062,0.285){\footnotesize \emph{MFG equations} }
\usefont{T1}{ptm}{m}{n}
\rput(11.883282,-0.295){\footnotesize \emph{2N equations}}
\end{pspicture}
}
\end{center}

This paper also opens the door to the important challenge of numerical resolution of Hamilton Jacobi equations in very high dimension either to solve the equations presented in this paper or as a numerical approximation to solutions of the Hamilton-Jacobi equations in infinite dimension presented in \cite{MFG4}.

\bibliographystyle{plain}
\nocite{*}
\bibliography{biblio}

\begin{thebibliography}{10}

\bibitem{achdou2010capuzzo}
Y.~Achdou, F.~Camilli, and I.~Capuzzo-Dolcetta.
\newblock Mean field games: numerical methods for the planning problem.
\newblock {\em Preprint}, 2010.

\bibitem{achdou2010mean}
Y.~Achdou and I.~Capuzzo-Dolcetta.
\newblock Mean field games: Numerical methods.
\newblock {\em SIAM Journal on Numerical Analysis}, 48(3):1136--1162, 2010.

\bibitem{cardaliaguet2010}
P.~Cardaliaguet.
\newblock Notes on mean field games.
\newblock 2010.

\bibitem{evans}
L.C. Evans.
\newblock {Partial Differential Equations (Graduate Studies in Mathematics,
  Vol. 19)}.
\newblock 2009.

\bibitem{gomes2010discrete}
D.A. Gomes, J.~Mohr, and R.R. Souza.
\newblock Discrete time, finite state space mean field games.
\newblock {\em Journal de Mathématiques Pures et Appliquées}, 93(3):308--328,
  2010.

\bibitem{gueant2008mean}
O.~Guéant.
\newblock {\em Mean field games and applications to economics}.
\newblock PhD thesis, PhD thesis. Université Paris-Dauphine, 2009.

\bibitem{gueant2009reference}
O.~Guéant.
\newblock A reference case for mean field games models.
\newblock {\em Journal de mathématiques pures et appliquées}, 92(3):276--294,
  2009.

\bibitem{gueant2011}
O.~Guéant.
\newblock Mean field games equations with quadratic hamiltonian: a specific
  approach.
\newblock {\em working paper}, 2011.

\bibitem{ParisPrinceton}
O.~Guéant, J.M. Lasry, and P.L. Lions.
\newblock Mean field games and applications.
\newblock In {\em Paris Princeton Lectures on Mathematical Finance}, 2010.

\bibitem{lachapelle2010computation}
A.~Lachapelle, J.~Salomon, and G.~Turinici.
\newblock Computation of mean field equilibria in economics.
\newblock {\em Mathematical Models and Methods in Applied Sciences}, 20(4):567,
  2010.

\bibitem{MFG1}
J.-M. Lasry and P.-L. Lions.
\newblock Jeux à  champ moyen i. le cas stationnaire.
\newblock {\em C. R. Acad. Sci. Paris}, 343(9), 2006.

\bibitem{MFG2}
J.-M. Lasry and P.-L. Lions.
\newblock Jeux à  champ moyen ii. horizon fini et contrôle optimal.
\newblock {\em C. R. Acad. Sci. Paris}, 343(10), 2006.

\bibitem{MFG3}
J.-M. Lasry and P.-L. Lions.
\newblock Mean field games.
\newblock {\em Japanese Journal of Mathematics}, 2(1), Mar. 2007.

\bibitem{MFG4}
P.-L. Lions.
\newblock Cours au collège de france : Théorie des jeux à champs moyens.
\newblock {\em
  http://www.college-de-france.fr/default/EN/all/equ\_der/audio\_video.jsp}.

\end{thebibliography}

\end{document}